\documentclass[12pt]{article}

\usepackage{amsmath,amsfonts,amssymb,latexsym,amscd}

 \input{prepictex.tex}
 \input{pictex.tex}
 \input{postpictex.tex}

\begin{document}

\newtheorem{theorem}{Theorem}[section]
\newtheorem{corollary}[theorem]{Corollary}
\newtheorem{lemma}[theorem]{Lemma}
\newtheorem{proposition}[theorem]{Proposition}
\newtheorem{conjecture}[theorem]{Conjecture}
\newtheorem{commento}[theorem]{Comment}
\newtheorem{definition}[theorem]{Definition}
\newtheorem{problem}[theorem]{Problem}
\newtheorem{remark}[theorem]{Remark}
\newtheorem{remarks}[theorem]{Remarks}
\newtheorem{example}[theorem]{Example}

\newcommand{\Nb}{{\mathbb{N}}}
\newcommand{\Rb}{{\mathbb{R}}}
\newcommand{\Tb}{{\mathbb{T}}}
\newcommand{\Zb}{{\mathbb{Z}}}
\newcommand{\Cb}{{\mathbb{C}}}

\newcommand{\Ef}{\mathfrak E}
\newcommand{\Gf}{\mathfrak G}
\newcommand{\iGf}{\mathfrak I\mathfrak G}
\newcommand{\Hf}{\mathfrak H}
\newcommand{\Kf}{\mathfrak K}
\newcommand{\Lf}{\mathfrak L}
\newcommand{\Af}{\mathfrak A}
\newcommand{\Bf}{\mathfrak B}

\def\A{{\mathcal A}}
\def\B{{\mathcal B}}
\def\C{{\mathcal C}}
\def\D{{\mathcal D}}
\def\F{{\mathcal F}}
\def\G{{\mathcal G}}
\def\H{{\mathcal H}}
\def\J{{\mathcal J}}
\def\K{{\mathcal K}}
\def\LL{{\mathcal L}}
\def\N{{\mathcal N}}
\def\M{{\mathcal M}}
\def\N{{\mathcal N}}
\def\O{{\mathcal O}}
\def\P{{\mathcal P}}
\def\SS{{\mathcal S}}
\def\T{{\mathcal T}}
\def\U{{\mathcal U}}
\def\W{{\mathcal W}}
\def\Z{{\mathcal Z}}

\def\span{\operatorname{span}}
\def\Ad{\operatorname{Ad}}
\def\ad{\operatorname{Ad}}
\def\tr{\operatorname{tr}}
\def\id{\operatorname{id}}
\def\en{\operatorname{End}}
\def\aut{\operatorname{Aut}}
\def\out{\operatorname{Out}}
\def\inn{\operatorname{Inn}}
\def\per{\operatorname{Per}(X_n)}
\def\la{\langle}
\def\ra{\rangle}

\def\j{j_\infty}
\def\f{f_\infty}
\def\g{g_\infty}
\def\a{a_\infty}

\title{The Weyl Group of the Cuntz Algebra}

\author{Roberto Conti\footnote{This research was supported through 
the programme "Research in Pairs" by the Mathematisches Forschungsinstitut Oberwolfach in 2011.}, 
Jeong Hee Hong\footnote{This research was supported by Basic Science Research Program through 
the National Research Foundation of Korea (NRF) funded by the Ministry of Education, Science and 
Technology (Grant No. 2012R1A1A2039991).} and
Wojciech Szyma{\'n}ski*\footnote{Partially supported by the FNU Rammebevilling
`Operator algebras and applications',  the NordForsk Research 
Network `Operator algebra and dynamics' (grant \#11580), and the FNU Forskningsprojekt 
`Structure and Symmetry'.}}

\date{}

\maketitle

\renewcommand{\sectionmark}[1]{}

\vspace{7mm}
\begin{abstract}
The Weyl group of the Cuntz algebra $\O_n$ is investigated.
This is (isomorphic to) the group of polynomial automorphisms $\lambda_u$ of 
$\O_n$, namely those induced by unitaries $u$ that can be written as 
finite sums of words in the canonical generating isometries $S_i$ and their adjoints.
A necessary and sufficient algorithmic combinatorial condition
is found for deciding when a polynomial endomorphism $\lambda_u$ restricts 
to an automorphism of the canonical diagonal MASA.
Some steps towards a general criterion for invertibility of $\lambda_u$ on the whole of $\O_n$ are also taken. 
A condition for verifying invertibility of a certain subclass of polynomial endomorphisms is given.
First examples of polynomial automorphisms of $\O_n$ not inner 
related to those induced by unitaries from the core UHF subalgebra are exhibited, for every $n\geq 2$.
In particular, the image of the Weyl group in the outer automorphism 
group of $\O_n$ is strictly larger than the image of the restricted Weyl group
analyzed in previous papers. Results about the action of the Weyl group on the spectrum of the 
diagonal are also included.
\end{abstract}

\vfill\noindent {\bf MSC 2010}: Primary 46L40, 46L55, Secondary 37B10

\vspace{3mm}

\newpage


\section{Introduction}

Consider a finite alphabet $\{1,2,\ldots,n\}$ with $n\geq2$ letters, and let ${\mathcal W}$ be the set of 
finite words on this alphabet. We say that two words are orthogonal if one is not the initial subword 
of the other. Let $\Sigma$ be the collection of finite subsets of ${\mathcal W}$ consisting of mutually 
orthogonal words. We consider 
the set of $n$ words $\{\alpha1,\alpha2,\ldots,\alpha n\}$ (beginning with the same subword $\alpha$ and 
ending with all the distinct letters of the alphabet) equivalent to the single word $\alpha$, and this 
extends to an equivalence relation on $\Sigma$. The set of equivalence classes is denoted $\widetilde{\Sigma}$. 
We fix $U\in\Sigma$, comprised of two ordered subsets: $\{\alpha_1,\ldots,\alpha_r\}$ and 
$\{\beta_1,\ldots,\beta_r\}$, with the 
property that both $\{\alpha_1,\ldots,\alpha_r\}$ and $\{\beta_1,\ldots,\beta_r\}$ are equivalent to the 
empty word. Such a $U$ determines recursively a sequence of transformations $T_k:\widetilde{\Sigma}\to
\widetilde{\Sigma}$ such that: 

\vspace{2mm}
 if $\gamma=\alpha_j\mu$ for some $j$ then $T_1(\gamma)=\beta_j\mu$, and 

\vspace{2mm} 
 if $T_{k-1}(\gamma)=\nu\alpha_j\mu$ for some $j$ and a word $\nu$ of length $k-1$ then 
$T_k(\gamma)=\nu\beta_j\mu$. 

\vspace{2mm}\noindent
Thus each transformation $T_k$ is determined by a certain Turing machine, 
\cite{E}, and hence it is computable for any finite set of inputs. 
We are interested in the following {\bf stabilization problem} of the recursive process $T_k\circ T_{k-1}\circ
\ldots\circ T_1:\widetilde{\Sigma}\to\widetilde{\Sigma}$. For what $U$ 
it holds that for each $A\in\widetilde{\Sigma}$ there exists an $m$ such that for all $k\geq m$ we have 
$T_k\circ\ldots\circ T_m\circ\ldots\circ T_1=T_m\circ\ldots\circ T_1$? 
We provide a surprisingly simple complete solution to this (suitably reformulated in more 
algebraic terms) stabilization problem in Theorem \ref{diagonaltheorem}, below. 

We can reformulate the above described combinatorial setup in topological terms, as follows. Let $X_n$ be the 
space of all (one-sided) infinite words. Then $X_n$ is a Cantor set with the product topology and elements 
of $\widetilde{\Sigma}$ are in bijective correspondence with its clopen subsets. Our stabilization problem 
is then equivalent to injectivity of a certain continuous map $\psi_U:X_n\to X_n$ determined naturally by $U$. 
Then, by the Gelfand duality, this problem is equivalent to surjectivity of a unital, injective $*$-homomorphism 
$\tilde{\psi}_U:C(X_n)\to C(X_n)$, dual to $\psi_U$. That is, we ask if $\psi_U$ is a homeomorphism of $X_n$ 
or, equivalently, if $\tilde{\psi}_U$ is an automorphism of $C(X_n)$. 

Somewhat paradoxically, it is most natural to view this problem in the context of much larger and  
noncommutative Cuntz algebras $\O_n$, \cite{Cun1}. These are $C^*$-algebras 
generated by $n$ isometries $S_1,\ldots,S_n$ of a Hilbert space with ranges adding up to the identity. 
In our setting, the element $U$ gives rise to a unitary $u=\sum_{j=1}^r S_{\alpha_j} S_{\beta_j}^*$ in 
$\O_n$, which in turn leads to a necessarily injective, unital $*$-endomorphism 
$\lambda_u$ of $\O_n$ such that $\lambda_u(S_j)=uS_j$ for all $j=1,\ldots,n$. The $C^*$-subalgebra 
$\D_n$ of $\O_n$ generated by ranges of all finite products of $S_1,\ldots,S_n$ is maximal abelian in $\O_n$ 
and naturally isomorphic to $C(X_n)$. The restriction of endomorphism $\lambda_u$ to $\D_n$ coincides 
with $\tilde{\psi}_U$. Thus, our combinatorial stabilization problem is equivalent to the {\bf problem of 
surjectivity} of $\lambda_u|_{\D_n}$. The question of surjectivity of $\lambda_u$ itself is very 
interesting as well and closely related to the so called Weyl group of the Cuntz algebra. This last problem 
appears very difficult and an algorithm for deciding surjectivity of an arbitrary $\lambda_u$ has not been found 
yet, although we make some headway towards its solution, below.  

\medskip
The present paper is a  continuation of our investigations of the subgroup $\aut(\O_n,\D_n)$  
of automorphisms of $\O_n$ which globally preserve the canonical diagonal MASA $\D_n$, and 
of related endomorphisms of $\O_n$, \cite{CS,CKS,CRS,HSS,Co,CHS1}. As shown in \cite{Cun2}, 
the quotient of $\aut(\O_n,\D_n)$ by its normal subgroup $\aut_{\D_n}(\O_n)$, consisting of those 
automorphisms which fix $\D_n$ point-wise, is discrete. Since $\aut_{\D_n}(\O_n)$
is a maximal abelian subgroup of $\aut(\O_n)$, \cite{Cun2}, it is natural to call this quotient  
the {\em Weyl group} of $\O_n$. The Weyl group contains a natural interesting subgroup 
corresponding to those automorphisms which also globally preserve the core UHF-subalgebra 
$\F_n$ of $\O_n$, called the {\em restricted Weyl group} of $\O_n$. It was shown in \cite{CHS1} 
that the image of the restricted Weyl group in the outer automorphism group of $\O_n$ 
can be embedded into the quotient of the automorphism group of the full two-sided 
$n$-shift by its center, and this embedding is surjective whenever $n$ is prime. In the present 
article, we focus our attention on the (full) Weyl group. It was shown in \cite{CS} that the Weyl group is 
isomorphic with the group of those automorphisms $\lambda_u\in\aut(\O_n)$ whose corresponding unitaries 
$u$ may be written as  a sum of words in $\{S_i,S_j^*\}$. (The collection of all such unitaries in $\O_n$  
is denoted $\SS_n$.) The structure of the Weyl group is highly complicated. For example, it contains the 
Thompson $F$ group in its intersection with $\inn(\O_n)$, \cite{N}. Our main objective  here
is investigation of the structure of the Weyl group of $\O_n$, its action on the diagonal MASA, 
and determining which unitaries $u\in\SS_n$ give rise to automorphisms.

The present paper is organized as follows. In section 2, we set up notation and review some basic facts 
on Cuntz algebras and their endomorphisms. In section 3, we study the restriction of an endomorphism 
$\lambda_u$, $u\in\SS_n$, to the diagonal $\D_n$. We give an 
algorithmic criterion for $\lambda_u|_{\D_n}$ to be an automorphism of $\D_n$, Theorem 
\ref{diagonaltheorem}. Its proof is combinatorial and involves equivalence of surjectivity of 
$\lambda_u|_{\D_n}$ with the stabilization problem mentioned above. In section 4, we investigate the 
problem when $\lambda_u$ is an automorphism of the entire $\O_n$. In Proposition \ref{invcriterion}, 
we present a combinatorial procedure for deciding this question for a certain large class of unitaries 
$u\in\SS_n$. In section 5, 
we exhibit endomorphisms $\lambda_u$, $u\in\SS_n$, which are not inner related to the ones of 
the form $\lambda_w$ with $w$ a unitary in the core UHF-subalgebra $\F_n$. In particular, we show 
with concrete examples that the image in $\out(\O_n)$ of the Weyl group is strictly larger then 
the image of the restricted Weyl group, Theorem \ref{outermain} and Corollary \ref{fullrestricted}. 
Finally, in section 6, we look at the action induced by $\lambda_u$ on the space $X_n$, the spectrum of 
the diagonal $\D_n$. We characterize homeomorphisms 
of $X_n$ corresponding to automorphisms $\Ad(u)$, $u\in\SS_n$, and describe the fixed points in $X_n$ 
for some exotic automorphisms $\lambda_u$. 


\section{Notation and preliminaries}

If $n$ is an integer greater than 1, then the Cuntz algebra $\O_n$ is a unital, simple,
purely infinite $C^*$-algebra generated by $n$ isometries $S_1, \ldots, S_n$ satisfying
$\sum_{i=1}^n S_i S_i^* = 1$, \cite{Cun1}.
We denote by $W_n^k$ the set of $k$-tuples $\mu = (\mu_1,\ldots,\mu_k)$
with $\mu_m \in \{1,\ldots,n\}$, and by $W_n$ the union $\cup_{k=0}^\infty W_n^k$,
where $W_n^0 = \{0\}$. We call elements of $W_n$ multi-indices.
If $\mu \in W_n^k$ then $|\mu| = k$ is the length of $\mu$. For $\mu,\nu\in W_n$ 
we write $\mu\prec\nu$ if $\mu$ is an initial subword of $\nu$.  If $\mu\in W_n^k$, 
$\nu\in W_n^m$ and $\mu\prec\nu$, then 
we denote by $\nu-\mu$ the word in $W_n^{m-k}$ obtained from $\nu$ by removing 
its initial segment $\mu$. Also, if $\mu\in W_n^k$ then we denote by $s(\mu)$ its first letter, 
and by $\tilde{\mu}$ the word in $W_n^{k-1}$ obtained from $\mu$ 
by removing $s(\mu)$. We denote by $\mu\wedge\nu$ the collection of all non-empty 
words $\eta$ such that both $\eta\prec\mu$ and $\eta\prec\nu$. 
If $\mu = (\mu_1,\ldots,\mu_k) \in W_n$ then $S_\mu = S_{\mu_1} \ldots S_{\mu_k}$
($S_0 = 1$ by convention) is an isometry with range projection $P_\mu=S_\mu S_\mu^*$. 
Every word in $\{S_i, S_i^* \ | \ i = 1,\ldots,n\}$ can be uniquely expressed as
$S_\mu S_\nu^*$, for $\mu, \nu \in W_n$ \cite[Lemma 1.3]{Cun1}.

We denote by $\F_n^k$ the $C^*$-subalgebra of $\O_n$ spanned by all words of the form
$S_\mu S_\nu^*$, $\mu, \nu \in W_n^k$, which is isomorphic to the
matrix algebra $M_{n^k}({\mathbb C})$. The norm closure $\F_n$ of
$\cup_{k=0}^\infty \F_n^k$ is the UHF-algebra of type $n^\infty$,
called the core UHF-subalgebra of $\O_n$, \cite{Cun1}. We denote by $\tau$ the 
unique normalized trace on $\F_n$. The core UHF-subalgebra $\F_n$ is the fixed-point 
algebra for the gauge action $\gamma:U(1)\to\aut(\O_n)$, such that $\gamma_z(S_j)=zS_j$ for 
$z\in U(1)$ and $j=1,\ldots,n$. We denote by $E$ the faithful conditional expectation from 
$\O_n$ onto $\F_n$ given by averaging with respect to the normalized Haar measure: 
$$ E(x) = \int_{z\in U(1)}\gamma_z(x)dz. $$
For an integer $m\in\Zb$ we denote $\O_n^{(m)}:=\{x\in\O_n : \gamma_z(x)=z^m x, \, 
\forall z\in U(1)\}$, a spectral subspace for $\gamma$. Then $\O_n^{(0)}=\F_n$ and 
for each positive integer $m$ and each 
$\alpha\in W_n^m$ we have $\O_n^{(m)}=\F_nS_\alpha$ and $\O_n^{(-m)}=S_\alpha^*\F_n$. 

The $C^*$-subalgebra of $\O_n$ generated by projections $P_\mu$, $\mu\in W_n$, is a 
MASA (maximal abelian subalgebra) in $\O_n$. We call it the {\em diagonal} and denote $\D_n$. 
Every projection in $\D_n$ of the form $P_\alpha$ for some $\alpha\in W_n$ 
will be called {\em standard}. The spectrum of $\D_n$ is naturally identified with 
$X_n$ --- the full one-sided $n$-shift space. For $d\in\D_n$ we denote by $M_d$ a map 
$M_d:\D_n\to \D_n$ such that $M_d(x)=dx$. 

As shown by Cuntz in \cite{Cun2}, there exists the following bijective correspondence
between unitaries in $\O_n$ (whose collection is denoted $\U(\O_n)$) and unital $*$-endomorphisms 
of $\O_n$ (whose collection we denote
$\en(\O_n)$). A unitary $u\in\U(\O_n)$ determines an endomorphism $\lambda_u$ by
$$ \lambda_u(S_i) = u S_i, \;\;\; i=1,\ldots, n. $$
Conversely, if $\rho :\O_n\rightarrow \O_n$ is an endomorphism, then
$\sum_{i=1}^n\rho(S_i)S_i^*=u$ gives a unitary $u\in\O_n$
such that $\rho=\lambda_u$. 
Composition of endomorphisms corresponds to a `convolution'
multiplication of unitaries as follows:
\begin{equation}\label{convolution}
\lambda_u \circ \lambda_w = \lambda_{\lambda_u(w)u}.
\end{equation}
If $A$ is either a unital $C^*$-subalgebra of $\O_n$ or a subset of $\U(\O_n)$, then we denote 
$\lambda(A)=\{\lambda_u\in\en(\O_n) : u\text{ unitary in }A\}$ and $\lambda(A)^{-1}=
\{\lambda_u\in\aut(\O_n) : u\text{ unitary in }A\}$. 

We denote by $\varphi$ the canonical shift:
$$ \varphi(x)=\sum_{i=1}^n S_ixS_i^*, \;\;\; x\in\O_n. $$
If we take $u=\sum_{i, j=1}^n S_iS_jS_i^*S_j^*$ then $\varphi=\lambda_u$. 
For all $u\in\U(\O_n)$ we have $\Ad(u)=\lambda_{u\varphi(u^*)}$. 
It is well-known that $\varphi$ leaves $\D_n$ globally invariant. We denote by 
$\phi$ the standard left inverse of $\varphi$, defined as 
$$ \phi(x)=\frac{1}{n}\sum_{i=1}^n S_i^* x S_i, \;\;\; x\in\O_n. $$ 
If $u\in\U(\O_n)$ then for each positive integer $k$ we denote
\begin{equation}\label{uk}
u_k = u \varphi(u) \cdots \varphi^{k-1}(u).
\end{equation}
Here $\varphi^0=\id$, and we agree that $u_k^*$ stands for $(u_k)^*$. If
$\alpha$ and $\beta$ are multi-indices of length $k$ and $m$, respectively, then
$\lambda_u(S_\alpha S_\beta^*)=u_kS_\alpha S_\beta^*u_m^*$. This is established through
a repeated application of the identity $S_i x = \varphi(x)S_i$, valid for all
$i=1,\ldots,n$ and $x \in \O_n$. 

We often consider elements of $\O_n$ of the form $w=\sum_{(\alpha,\beta)\in\J}c_{\alpha,\beta}
S_\alpha S_\beta^*$, where $\J$ is a finite collection  of pairs $(\alpha,\beta)$ of words $\alpha,\beta\in 
W_n$ and $c_{\alpha,\beta}\in\Cb$. We denote $\J_1=\{\alpha:\exists (\alpha,\beta)\in\J\}$ and 
$\J_2=\{\beta:\exists (\alpha,\beta)\in\J\}$. 
Of course, such a presentation (if it exists) is not unique, but once it is 
chosen then we associate with it two integers: $\ell=\ell(\J)=\max\{|\alpha|:(\alpha,\beta)\in\J\}$ 
and $\ell'=\ell'(\J)=\max\{|\alpha|,|\beta|:(\alpha,\beta)\in\J\}$. 
Note that if $w\in\F_n$ then $w\in\F_n^\ell$.  We have 
$\varphi(w)=\sum_{(\mu,\nu)\in\varphi(\J)}c_{\mu,\nu}S_\mu S_\nu^*$, where 
$\varphi(\J)=\{((i,\alpha),(\beta,i)) : i\in W_n^1,\; (\alpha,\beta)\in\J\}$ and 
$c_{(i,\alpha),(\beta,i)}=c_{\alpha,\beta}$. Then $\ell(\varphi(\J))=\ell(\J)+1$ and 
$\ell'(\varphi(\J))=\ell'(\J)+1$. In particular, we consider the group $\SS_n$ 
 of those unitaries in $\O_n$ which can be written
as finite sums of words, i.e. in the form $u=\sum_{(\alpha,\beta)\in\J}S_\alpha S_\beta^*$. Note that 
such a sum is a unitary if and only if $\sum_{\alpha\in\J_1} P_\alpha = 1 = \sum_{\beta\in\J_2} P_\beta$. 
We also write $\P_n=\SS_n\cap\F_n$ and $\P_n^k=\SS_n\cap\F_n^k$ for the 
subgroups of $\SS_n$ consisting of permutative unitaries. 

For algebras $A\subseteq B$ we denote by $\N_B(A)=\{u\in\U(B):uAu^*=A\}$ the normalizer
of $A$ in $B$ and by $A' \cap B=\{b \in B: (\forall a \in A) \; ab=ba\}$ the 
relative commutant of $A$ in $B$. We also denote by $\aut(B,A)$ the collection of all those 
automorphisms $\alpha$ of $B$ such that $\alpha(A)=A$, and by $\aut_A(B)$ those 
automorphisms of $B$ which fix $A$ point-wise. 

$\aut_{\D_n}(\O_n)$ is a normal subgroup of $\aut(\O_n,\D_n)$, and the corresponding quotient 
is called the {\em Weyl group} of $\O_n$. It was shown in \cite{Cun2} that the Weyl group is discrete, 
and more recently in \cite{CS} that it is isomorphic to $\lambda(\SS_n)^{-1}$. The quotient of 
$\aut(\O_n,\D_n)\cap\aut(\O_n,\F_n)$ by $\aut_{\D_n}(\O_n)$ is called the {\em restricted 
Weyl group} of $\O_n$. It is isomorphic to $\lambda(\P_n)^{-1}$, \cite{CS}. The image of $\lambda(\SS_n)^{-1}$ 
in $\out(\O_n)$ is called the {\em outer Weyl group} of $\O_n$ and such image of $\lambda(\P_n)^{-1}$ 
is called the {\em restricted outer Weyl group} of $\O_n$. As shown in \cite[Theorem 3.7]{CHS2}, 
the outer Weyl group is just the quotient of $\lambda(\SS_n)^{-1}$ by $\{\Ad(u):u\in\SS_n\}$. 
Likewise, the restricted outer Weyl group is the quotient of $\lambda(\P_n)^{-1}$ by $\{\Ad(w):w\in\P_n\}$. 


\section{The automorphisms of the diagonal}

In this section, we give an algorithmic crierion for deciding if the restriction to $\D_n$ 
of an endomorphisms $\lambda_u$, $u\in\SS_n$, gives rise to an automorphism of the diagonal $\D_n$. 

\begin{lemma}\label{push}
Let $u\in\SS_n$ be such that $u=\sum_{(\alpha,\beta)\in\J}S_\alpha S_\beta^*$, and let  
$\ell=\ell(\J)$. Then $\lambda_u(\D_n^k)\subseteq\D_n^{k\ell}$ for all $k\in\Nb$. 
\end{lemma}
{\em Proof.} 
We proceed by induction on $k$. For $k=1$ and $i\in W_n^1$ we have
$$ \lambda_u(P_i)=uP_i u^*=\sum_{(\alpha,\beta),\;(\alpha',\beta')\in\J}S_\alpha S_\beta^* 
P_i S_{\beta'} S_{\alpha'}^*= \sum_{(\alpha,\beta)\in\J, \; s(\beta)=i}S_\alpha S_\alpha^* $$
and thus $\lambda_u(\D_n^1)\subseteq\D_n^\ell$. For the inductive step, suppose that 
$\lambda_u(\D_n^k)\subseteq\D_n^{k\ell}$. Then 
$$ \lambda_u(\D_n^{k+1})=\lambda_u(\D_n^1\varphi(\D_n^k))
=\lambda_u(\D_n^1) (\Ad(u) \varphi\lambda_u)(\D_n^k) 
\subseteq \D_n^\ell (\Ad(u) \varphi)(\D_n^{k\ell}) 
\subseteq \D_n^\ell \D_n^{(k+1)\ell} $$
and thus $\lambda_u(\D_n^{k+1})\subseteq\D_n^{(k+1)\ell}$. 
\hfill$\Box$

\begin{proposition}\label{diagonalinvert}
Let $u\in\SS_n$. Then the following hold. 
\begin{enumerate}
\item $\lambda_u|_{\D_n}$ is an automorphism of $\D_n$ if and only if for each $\alpha\in W_n$ 
the sequence $\{u_k^* P_\alpha u_k\}$ eventually stabilizes. 
\item $\lambda_u$ is an automorphism of $\O_n$ if and only if: 
\begin{enumerate}
\item $\lambda_u|_{\D_n}$ is an automorphism of $\D_n$, and 
\item there exists a $w\in\SS_n$ such that $\lambda_w|_{\D_n}=(\lambda_u|_{\D_n})^{-1}$. 
\end{enumerate}
\end{enumerate}
\end{proposition}
{\em Proof.} 
Ad 1. This is well-known, \cite{Cun2}. Indeed, the sequence $\{u_k^* P_\alpha u_k\}$ eventually stabilizes 
if and only if $P_\alpha$ belongs to the range of $\lambda_u$ (and then $\lambda_u(\lim 
u_k^* P_\alpha u_k)=P_\alpha$). Thus, condition 1. is equivalent to $\lambda_u(\D_n)=\D_n$, i.e. to 
$\lambda_u|_{\D_n}$ being an automorphism of $\D_n$. 

Ad 2. If $\lambda_u$ is automorphism of $\O_n$, then $\lambda_u(\D_n)\subseteq\D_n$ since $u\in
\N_{\O_n}(\D_n)$. Thus $\lambda_u(\D_n)=\D_n$, since $\D_n$ is a MASA in $\O_n$. Also,  
there exists $w\in\SS_n$ such that $\lambda_u^{-1}=\lambda_w$, \cite{S,CS,M}. This gives one 
implication of part 2. For the reversed implication, suppose that (a) and (b) hold. Then $\lambda_u
\lambda_w|_{\D_n}=\id$. Thus $\lambda_u\lambda_w$ is an automorphism of $\O_n$ by 
\cite[Proposition 3.2]{Co}. Consequently, $\lambda_u$ being surjective is automorphism of $\O_n$. 
\hfill$\Box$

\begin{example}\label{examples123}
\rm {\bf (a)} If $u=S_1S_1S_1^*+S_1S_2S_1^*S_2^*+S_2S_2^*S_2^*\in\SS_2$  then 
$\lambda_u|_{\D_2}$ is not surjective. Indeed, a straightforward calculation shows that 
$\lambda_u(\D_2)P_2=\Cb P_2$. 

\vspace{2mm}\noindent
{\bf(b)} If $u=S_2S_1S_1^*+S_2S_2S_1^*S_2^*+S_1S_2^*S_2^*\in\SS_2$ then 
$\lambda_u|_{\D_2}$ is not surjective. Indeed, projection $P_{12}$ does not satisfy (Condition 
1) of Proposition \ref{diagonalinvert}. 

\vspace{2mm}\noindent
{\bf (c)} If $u=u^*=S_1S_2^*S_2^*+P_{21}+S_2S_2S_1^*\in\SS_2$ then 
$\lambda_u|_{\D_2}$ is not surjective. Indeed, projection $P_{11}$ does not satisfy (Condition 
1) of Proposition \ref{diagonalinvert}.
\end{example}

Our next result shows that in order to verify (Condition 1) in Proposition \ref{diagonalinvert}
it is enough to check it only for finitely many projections. Before that, we note the following. 
Let $u\in\SS_n$ be such that $u=\sum_{(\alpha,\beta)\in\J}S_\alpha S_\beta^*$. Then for 
each word $\mu\in W_n$ and for each $(\alpha,\beta)\in\J$ we have 
\begin{equation}\label{adaction}
\Ad(u)(P_{\beta\mu}) = P_{\alpha\mu}. 
\end{equation}
In particular, $\Ad(P_\beta)=P_\alpha$.

\begin{lemma}\label{localauto}
Let $u\in\SS_n$ be such that $u=\sum_{(\alpha,\beta)\in\J}S_\alpha S_\beta^*$, and let  
$\ell' = \ell'(\J)$. Then $\lambda_u|_{\D_n}$ is an automorphism of $\D_n$ if and only if for 
each $\gamma\in W_n^{\ell'}$ the sequence $\{u_k^* P_\gamma u_k\}$ eventually stabilizes. 
\end{lemma}
{\em Proof.}
For short, say a projection $Q \in \D_n$ is ``bad'' (relative to $u$) if the sequence 
$\{u_k^* Q u_k\}$ does not stabilize, and ``good'' otherwise. Also, let $r$ be the non-negative 
integer uniquely defined by requiring that all projections in $\D_n^r$
are good, but there is a bad projection in $\D_n^{r+1}$. Then at least one of the minimal projections 
in $\D_n^{r+1}$ is bad as well. We claim that $r+1 \leq \ell'$.

Reasoning by way of contradiction, suppose that $\ell' < r+1$ and let $p=P_{\gamma}$, 
$\gamma \in W_n^{r+1}$, be such a bad minimal projection in $\D_n^{r+1}$. Now, $u^* p u$ can be 
computed using equation (\ref{adaction}), with $u$ replaced by $u^*$, and hence it is still of the 
form $P_{\gamma_1}$ for some $\gamma_1 \in W_n$. In this process, by replacing the initial 
$\alpha$-segment of $\gamma$ with the corresponding $\beta$, the last $r+1 - \ell'$ digits will remain 
unaltered. Now, the assumption that $p$ is bad easily implies that the projection 
$P_\delta:=n\phi(u^* p u)$, obtained from $u^* p u$ by deleting the first digit of $\gamma_1$, is 
still bad. By assumption, one must have $|\delta| \geq r+1$, and hence $u^* p u \notin \D_n^{r+1}$. 
In other words, when computing $u^* p u$ we have replaced a word $\alpha$ in $\gamma$ with 
a longer word $\beta$. This implies that when in the next step we consider 
$ P_{\gamma_2}:=\varphi(u)^* u^* p u \varphi(u)$, 
the last $r+1 - \ell'$ digits of $\gamma_2$ will coincide again with those of $\gamma$.
Also, $(n\phi)^2(P_{\gamma_2})$ must be bad, i.e. $u_2^*pu_2=P_{\gamma_2} \notin \D_n^{r+2}$.
Repeating this argument, one can indeed show that $u_k^* p u_k = 
P_{\gamma_k}\notin \D_n^{r+k}$ for all $k=1,2,\ldots$,
and moreover the last $r+1 - \ell'$ digits of $\gamma_k$ coincide with those of $\gamma$ for any $k$.
All in all, this means that these last digits of $\gamma$ indeed play no role in the whole process and 
defining $\gamma'$ simply to be the multi-index obtained from $\gamma$ by deleting its last digit,
the very same argument would readily show that $P_{\gamma'}$ is still bad.
But then $P_{\gamma'} \in \D_n^r$, contradicting our assumption.

By the above, if there are bad projections at all, we can find at least one of them in 
$\D_n^{\ell'}$. As a sum of good projections is clearly good,
it is also clear that in that case there is always such a bad projection of the form $P_\gamma$, 
where $|\gamma| = \ell'$.
\hfill$\Box$

\medskip
All in all, for $u \in \SS_n$ one has 
\begin{equation}\label{finiteond}
\lambda_u(\D_n)=\D_n \;\Leftrightarrow\; \D_n^{\ell'} \subseteq \lambda_u(\D_n), 
\end{equation}
where $\ell'$ is as in the statement of Lemma \ref{localauto}.

\medskip
In view of Lemma \ref{localauto}, the process of determination if an endomorphism $\lambda_u|_{\D_n}$, 
$u\in\SS_n$, is an automorphism of the diagonal can be reduced to verification if a certain finite collection 
of projections is contained in its range. This is a very significant reduction but still it is not clear a priori if 
this process can be carried out in finately many steps even for a single projection! This question has 
a positive answer in the case of a permutative unitary $u\in\P_n$, as shown in \cite{S,CS}, but 
the present case is much more complicated. Now, we will describe a key construction of the present paper, 
producing a certain finite directed graph corresponding to a unitary $u\in\SS_n$. Non occurence of closed 
paths on the graph will turn out to be equivalent to $\lambda_u|_{\D_n}$ being automorphism of $\D_n$. 

\medskip
Given $u=\sum_{(\alpha,\beta)\in\J}S_\alpha S_\beta^*$ in $\SS_n$, we define a {\bf finite directed graph} 
$\Gamma_u$, whose vertices $\Gamma_u^0$ will be identified with certain subsets of $\J_1$. 
In order to construct the graph $\Gamma_u$, we proceed by induction. 

\vspace{2mm}\noindent
{\bf The initial step.} To begin with, we include in $\Gamma_u^0$ each singleton subset $\{\alpha\}$ of $\J_1$ 
and the empty set $\emptyset$. Now, given $(\alpha,\beta)\in\J$, one of the following three cases takes place:

\vspace{2mm}\noindent
(i) $\beta=(i)$ for some $i\in W_n^1$, 

\vspace{2mm}\noindent
(ii) $\beta=(i,\alpha',\mu)$ for some $i\in W_n^1$, $\alpha'\in\J_1$, and a word $\mu$ (possibly empty), 

\vspace{2mm}\noindent
(iii) $\beta=(i,\mu)$ for some $i\in W_n^1$ and a word $\mu$ which is an initial segment of at least 
two elements of $\J_1$, namely $\alpha_1',\ldots,\alpha_r'$. 

\vspace{2mm}
Depending on the case, we enlarge the graph $\Gamma_u$ as follows. In case (i), we add an edge from 
vertex $\{\alpha\}$ to vertex $\emptyset$ with label $i$. In case (ii), we add an edge from vertex 
$\{\alpha\}$ to vertex $\{\alpha'\}$ with label $i$. In case (iii), we add a vertex $A=\{\alpha_1',
\ldots,\alpha_r'\}$ and an edge from $\{\alpha\}$ to $A$ with label $i$. 

\vspace{2mm}\noindent
{\bf The inductive step.} Let $A\subseteq\J_1$ be a vertex added to $\Gamma_u^0$ in the preceding 
step, but $A\neq\emptyset$ and $A$ not a singleton set. For each $j\in W_n^1$ we proceed as follows. 
Let $B_k$, $k=1,\ldots,m$, be the collection of all those already constructed vertices of $\Gamma_u$ 
that there exists an $\alpha\in A$ and an edge from $\{\alpha\}$ to $B_k$ with label $j$. If 
$\bigcup_{k=1}^m B_k=\J_1$ then we add an edge from $A$ to $\emptyset$ with label $j$. 
If $\bigcup_{k=1}^m B_k\not=\J_1$ then we add a vertex $B=\bigcup_{k=1}^m B_k$ (if such 
a vertex does not exist already), and we add an edge from $A$ to $B$ with label $j$. 

\vspace{2mm}
Continuing inductively in the above described manner, we produce the desired graph $\Gamma_u$. 
This is a finite, directed, and labeled graph. 
Each vertex emits at most $n$ edges, carrying distinct labels from the set $W_n^1$. 
Any finite path on the graph $\Gamma_u$ may be uniquely identified with a pair $(A,\nu)$, where 
$A\in\Gamma_u^0$ is the initial vertex of the path and $\nu=(\nu_1,\nu_2,\ldots,\nu_k)$ is the word 
such that $\nu_j$ is the label of the $j^{\text th}$ edge entering this path. For such a path $(A,\nu)$, we 
denote its terminal vertex by $\nu(A)$. We will denote by $\Gamma_u^1$ 
the set of edges of the graph, by $\Gamma_u^k$ the set of paths of length $k$, and 
by $\Gamma_u^*$ the set of finite paths. $\Gamma_u^*(A)$  and 
$\Gamma_u^k(A)$, respectively, are the sets of finite paths and paths of length $k$ 
which begin at the vertex $A$. 

\begin{example}\label{graphwithloop}
\rm Let $u=S_1S_2^*S_2^*+S_2S_1S_1^*S_2^*+S_2S_2S_1^*$. Then the corresponding 
graph $\Gamma_u$ has five vertices and five edges, and looks as follows. In particular, there is a 
closed (directed) path on the graph. 

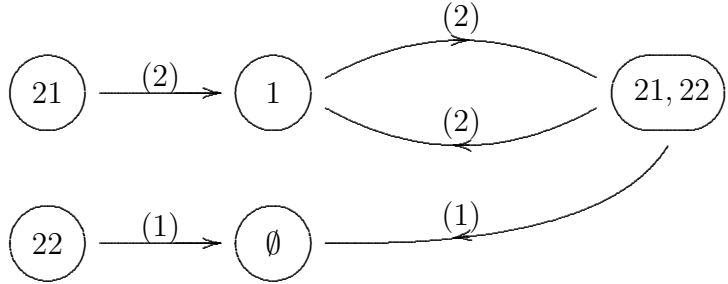
\[ \beginpicture
\setcoordinatesystem units <1cm,1cm>
\setplotarea x from -4 to 6, y from -2 to 2
\circulararc 360 degrees from -4 1 center at -3.5 1
\put {$21$} at -3.5 1
\circulararc 360 degrees from -1 1 center at -0.5 1
\put {$1$} at -0.5 1
\circulararc 360 degrees from -4 -1 center at -3.5 -1
\put {$22$} at -3.5 -1
\circulararc 360 degrees from -1 -1 center at -0.5 -1
\put {$\emptyset$} at -0.5 -1
\setlinear
\plot -2.8 1   -1.2 1 /
\plot -2.8 -1  -1.2 -1 /
\plot 4.5 1.5   5 1.5 /
\plot 4.5 0.5   5 0.5 /
\arrow <0.235cm> [0.2,0.6] from -1.8 1 to -1.2 1
\arrow <0.235cm> [0.2,0.6] from -1.8 -1 to -1.2 -1
\put {$(2)$} at -2 1.2
\put {$(1)$} at -2 -0.8
\setquadratic
\plot 0.2 1.2   2 1.7   3.8 1.2 /
\plot 0.2 0.8   2 0.3   3.8 0.8 /
\plot 0.2 -1   3.2 -0.7   4.75 0.3 /
\arrow <0.235cm> [0.2,0.6] from 1.9 1.7 to 2.1 1.7 
\arrow <0.235cm> [0.2,0.6] from 2.1 0.3 to 1.9 0.3 
\put {$(2)$} at 2 2
\put {$(2)$} at 2 0.6
\circulararc 180 degrees from 4.5 1.5 center at 4.5 1
\circulararc 180 degrees from 5 0.5 center at 5 1
\put {$21,22$} at 4.8 1
\arrow <0.235cm> [0.2,0.6] from 2.1 -0.9 to 1.9 -0.93
\put {$(1)$} at 2 -0.65
\endpicture \] 

\end{example}

\begin{example}\label{graphwithoutloop}
\rm Let $u=S_{12}S_{21}^*+S_{11}S_{221}^*+S_{21}S_{222}^*+S_{2222}S_{11}^*
+S_{2221}S_{122}^*+S_{221}S_{121}^*$. Then the corresponding 
graph $\Gamma_u$ looks as follows:

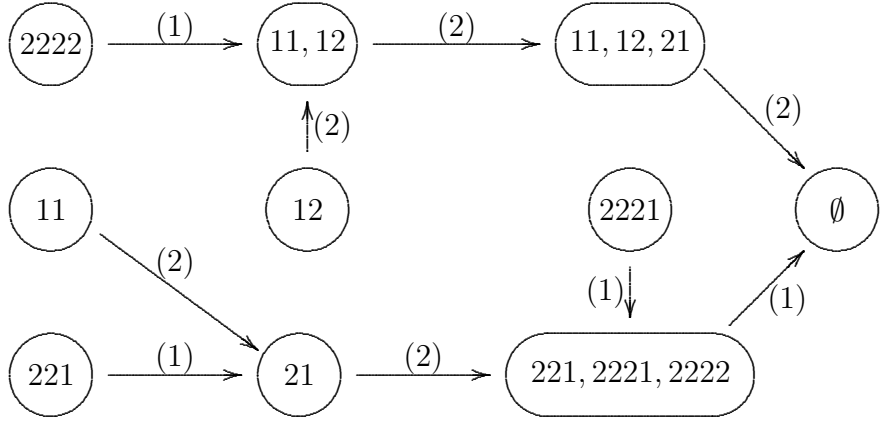
\[ \beginpicture
\setcoordinatesystem units <1.1cm,1.1cm>
\setplotarea x from 0 to 10.5, y from 0 to 5
\circulararc 360 degrees from  0.5 0 center at 0.5 0.5
\circulararc 360 degrees from 0.5 2 center at 0.5 2.5
\circulararc 360 degrees from  3.5 0 center at 3.5 0.5
\circulararc 180 degrees from 6.5 1 center at 6.5 0.5 
\circulararc 180 degrees from 8.5 0 center at 8.5 0.5
\circulararc 360 degrees from 0.5 5 center at 0.5 4.5
\circulararc 180 degrees from 3.5 5 center at 3.5 4.5
\circulararc 180 degrees from 3.7 4 center at 3.7 4.5 
\circulararc 360 degrees from 3.6 3 center at 3.6 2.5
\circulararc 180 degrees from 7.1 5 center at 7.1 4.5
\circulararc 180 degrees from 7.9 4 center at 7.9 4.5
\circulararc 360 degrees from 7.5 3 center at 7.5 2.5
\circulararc 360 degrees from 10 3 center at 10 2.5

\put {$221$} at 0.5 0.5
\put {$11$} at 0.5 2.5
\put {$21$} at 3.5 0.5
\put {$221,2221,2222$} at 7.5 0.5
\put {$2222$} at 0.5 4.5
\put {$11,12$} at 3.6 4.5
\put {$12$} at 3.6 2.5
\put {$11,12,21$} at 7.5 4.5
\put {$2221$} at 7.5 2.5
\put {$\emptyset$} at 10 2.5

\put {$(1)$} at 2 4.7
\put {$(1)$} at 2 0.7
\put {$(1)$} at 7.2 1.5
\put {$(2)$} at 5 0.7
\put {$(2)$} at 5.4 4.7
\put {$(2)$} at 3.9 3.5
\put {$(2)$} at 2 1.85
\put {$(2)$} at 9.35 3.7
\put {$(1)$} at 9.4 1.4

\setlinear
\plot 1.2 0.5   2.8 0.5 /
\plot 1.2 4.5   2.8 4.5 /
\plot 1.1 2.2   3 0.8 /
\plot 4.2 0.5   5.8 0.5 /
\plot 6.5 0   8.5 0 /
\plot 6.5 1   8.5 1 /
\plot 3.5 5   3.7 5 /
\plot 3.5 4   3.7 4 /
\plot 3.6 3.2   3.6 3.8 /
\plot 7.1 5   7.9 5 /
\plot 7.1 4   7.9 4 /
\plot 4.4 4.5   6.4 4.5 /   
\plot 7.5 1.8   7.5 1.2 /
\plot 8.4 4.2   9.6 3 /
\plot 8.7 1.1   9.6 2 /

\arrow <0.235cm> [0.2,0.6] from 2.6 0.5 to 2.8 0.5
\arrow <0.235cm> [0.2,0.6] from 2.6 4.5 to 2.8 4.5
\arrow <0.235cm> [0.2,0.6] from 2.81 0.94  to 3 0.8
\arrow <0.235cm> [0.2,0.6] from 5.6 0.5 to 5.8 0.5
\arrow <0.235cm> [0.2,0.6] from 3.6 3.6 to 3.6 3.8
\arrow <0.235cm> [0.2,0.6] from 6.2 4.5 to 6.4 4.5
\arrow <0.235cm> [0.2,0.6] from 7.5 1.4 to 7.5 1.2
\arrow <0.235cm> [0.2,0.6] from 9.4 3.2 to 9.6 3
\arrow <0.235cm> [0.2,0.6] from 9.4 1.8 to 9.6 2

\endpicture \]
 
\end{example}

For $\alpha\in\J_1$, we say that $\{\alpha\}$ is a 
{\em splitting vertex} if it emits an edge to a vertex $A\subseteq\J_1$ such that $A$ contains at 
least two elements. This happens when for $(\alpha,\beta)\in\J$ we have that $\beta$ is 
an initial subword of more than one $\alpha\in\J_1$. For example, $\alpha=1$ in Example 
\ref{graphwithloop} and $\alpha_1=12$, $\alpha_2=21$, $\alpha_3=2221$ and $\alpha_4=2222$ in 
Example \ref{graphwithoutloop} are all splitting vertices. 

The point of introducing graph $\Gamma_u$ is that it conveniently captures the essential 
features of the process of calculating $u_k^* P_\alpha u_k$, appearing in part 1 of 
Proposition \ref{diagonalinvert}. Indeed, for $A\in\Gamma_u^0$ denote $P_A:=
\sum_{\alpha\in A}P_\alpha$. Then we have 
\begin{equation}\label{u1pa}
\Ad(u^*)(P_A) = \sum_{(A,j)\in\Gamma_u^1(A)}\;\sum_{\alpha\in j(A)} P_{j\alpha}. 
\end{equation}
Clearly, if $\Ad(u^*)(P_\mu)=\sum_k P_{\nu_k}$, then for each $i\in W_n^1$ we have 
\begin{equation}\label{u2pa}
\Ad(\varphi(u^*))(P_{i\mu}) = \sum_k P_{i\nu_k}. 
\end{equation}
Combining (\ref{u1pa}) with (\ref{u2pa}) and proceeding by induction on $k$, we see that for any 
$A\in\Gamma_u^0$ and a non-negative integer $k$ we have 
\begin{equation}\label{ukpa}
\Ad(u^*_k)(P_A) = \sum_{(A,\nu)\in\Gamma_u^k(A)}\;\sum_{\alpha\in \nu(A)} P_{\nu\alpha}+
\sum_{m=1}^{k-1}\;\sum_{\substack{(A,\mu)\in\Gamma_u^{m}(A) \\ \mu(A)=\emptyset}}P_{\mu}. 
\end{equation} 

Now, we are ready to prove a theorem which gives an algorithmic (finite) procedure for 
determining if an endomorphism $\lambda_u$, $u\in\SS_n$, restricts to an automorphism 
of the diagonal $\D_n$. 

\begin{theorem}\label{diagonaltheorem}
Let $u\in\SS_n$ and let $\Gamma_u$ be the directed graph corresponding to $u$. Then 
$\lambda_u|_{\D_n}$ is an automorphism of $\D_n$ if and only if graph $\Gamma_u$ 
does not contain any closed (directed) paths. 
\end{theorem}
{\em Proof.}
Firstly, suppose that there is a closed path 
$$ A_1 \overset{(i_1)}{\longrightarrow} A_2 \overset{(i_2)}{\longrightarrow} 
\ldots \overset{(i_{r-1})}{\longrightarrow} A_r \overset{(i_r)}{\longrightarrow} A_1 $$
in the graph $\Gamma_u$. We denote $\nu=(i_1,i_2,\ldots,i_r)$ and $\nu^k=\nu\nu\cdots\nu$ ($k$-fold 
composition). With help of formula (\ref{ukpa}) we see that 
\begin{equation}\label{loop}
\Ad(u^*_{kr})(P_{A_1}) = P_{\nu^k}\varphi^{kr}(P_{A_1}) + 
\sum_{\mu\wedge\nu^k=\emptyset}P_\mu. 
\end{equation}
Given any $k<k'$ there exists a non-zero projection $q\in\D_n$ such that $q\leq P_{\nu^k}$ and 
$qP_{\nu^{k'}}=0$. But then formula (\ref{loop}) implies that $q\Ad(u^*_{kr})(P_{A_1})\neq 0$ 
while $q\Ad(u^*_{k'r})(P_{A_1})=0$. Thus the sequence $\{\Ad(u^*_m)(P_{A_1})\}$ never 
stabilizes and, consequently, projection $P_{A_1}$ does not belong to $\lambda_u(\D_n)$, Proposition 
\ref{diagonalinvert}. 

\vspace{2mm}
Conversely, suppose that graph $\Gamma_u$ does not contain any closed paths. By virtue of Lemma 
\ref{localauto}, it suffices to show that the sequence $\{\Ad(u^*_k)(P_\mu)\}$ eventually stabilizes 
for each $\mu\in W_n^{\ell'}$ with $\ell'=\ell'(\J)$. To this end, consider the following three cases. 

Firstly, we consider the case of $P_\alpha$, 
$\alpha\in\J_1$. Since $\Gamma_u$ is a finite graph without closed paths, there are only finitely 
many paths and each of them terminates at a sink. By construction, graph $\Gamma_u$ contains 
exactly one sink, namely vertex $\emptyset$. Thus formula (\ref{ukpa}) applied to $A=\{\alpha\}$ 
shows that for sufficiently large $k$ we have 
$$ \Ad(u^*_k)(P_\alpha) = \sum_{\substack{(A,\mu)\in\Gamma^*_u(\{\alpha\}) \\ 
\mu(\{\alpha\})=\emptyset}}P_\mu, $$
and thus the sequence $\{\Ad(u^*_k)(P_\alpha)\}$ eventually stabilizes. 

Secondly, we consider a word $\mu$ such that 
there exists an $\alpha\in\J_1$ with $\mu\prec\alpha$. Then $P_\mu=\sum P_{\alpha'}$, where the 
sum is over all such $\alpha'\in\J_1$ that $\mu\prec\alpha'$. In this case, the sequence 
$\{\Ad(u^*_k)(P_\mu)\}$ stabilizes by the preceding argument. 

Thirdly, we must consider the case with $\mu$ a word of length at most $\ell'$ for which there 
exists an $\alpha\in\J_1$ such that $\alpha\prec\mu$. Then write $\mu=\alpha\nu$. Let $(\{\alpha\},
\eta)$ be the maximal path beginning at $\{\alpha\}$ and such that each vertex on the path is a 
singleton subset of $\J_1$. Let $k$ be the length of this path. Using formula (\ref{ukpa}), we see that 
\begin{equation}\label{ukpmu}
\Ad(u^*_k)(P_\mu)=P_{\eta\alpha'\nu} 
\end{equation} 
for some $\alpha'\in\J_1$. Now, one of the following two cases happens: either $\{\alpha'\}$ emits 
an edge (with label $j$) to the sink $\emptyset$, or $\{\alpha'\}$ is a splitting vertex. In the former case, 
we have $\Ad(u^*_{k+1})(P_\mu)=P_{\eta j\nu}$, and the question of stabilization of the sequence 
corresponding to the word $\mu$ reduces to the same question for the sequence corresponding to 
the word $\nu$, which is strictly shorter then $\mu$. In the latter case, let $\{\alpha'\}$ emit an edge 
(with label $i$) to a vertex $A$. Then we have 
$\Ad(u^*)(\P_{\alpha'})=P_{\beta'}=\sum_{j=1}^m P_{i\alpha_j\nu_j}$, for 
some $\alpha_j\in\J_1$ and words $\nu_j$ such that each $\nu_j$ is 
strictly shorter then $\nu$. Taking into account formula (\ref{ukpmu}), we obtain 
$\Ad(u^*_{k+1})(P_\mu)=\sum_{j=1}^m P_{\eta i\alpha_j\nu_j}$. Thus, the question if 
the sequence $\{\Ad(u^*_k)(P_\mu)\}$ stabilizes (with $\mu=\alpha\nu$) reduces to the same question 
for all $\mu_j=\alpha_j\nu_j$, where $|\nu_j|<|\nu|$. 
Consequently, the claim follows for all words $\mu=\alpha\nu$, $\alpha\in\J_1$, by induction on $|\nu|$. 
\hfill$\Box$

\begin{remark}
\rm We note that for certain special classes of unitaries $u\in\SS_n$, different criteria for 
$\lambda_u|_{\D_n}\in\aut(\D_n)$ were given earlier in \cite{CS2}. 
\end{remark}


\section{The invertibility}

In this section, we consider the problem when $\lambda_u$, $u\in\SS_n$, is an automorphism of $\O_n$. 
Recall that $E:\O_n\to\F_n$ is the gauge invariant conditional expectation, and for a $\beta\in W_n^k$ the 
symbol $\tilde{\beta}$ denotes the word in $W_n^{k-1}$ obtained from $\beta$ by removing its first letter. 

\begin{lemma}\label{enz}
If $u\in\SS_n$ is arbitrary then there exists a $v\in\SS_n$ such that $E(w)\neq0$ 
for $w=vu\varphi(v^*)$. 
\end{lemma}
{\em Proof.}
Let $u=\sum_{(\alpha,\beta)\in\J}S_\alpha S_\beta^*$ and suppose that $E(u)=0$. If $v\in\SS_n$ 
then $vu\varphi(v^*)=\sum_{(\alpha,\beta)\in\J} vS_\alpha S_{\tilde{\beta}}^*v^*S_{\beta_1}^*$. 
Thus, it suffices to find a $v\in\SS_n$ such that for certain $(\alpha,\beta)\in\J$ we have 
$vS_\alpha S_{\tilde{\beta}}^*v^*\in\O_n^{(1)}$. Since $E(u)=0$, there exists $(\alpha,\beta)\in\J$ 
with $|\alpha|>|\beta|$. Now, one of the following two cases takes place: either $P_\alpha$ is 
orthogonal to $P_{\tilde{\beta}}$ or $\tilde{\beta}\prec\alpha$ and $\tilde{\beta}\neq\alpha$. 
In the former case, put $v=S_1^2S_\alpha^*+S_2S_{\tilde{\beta}}^*\;+$(other terms). In the latter, 
we have $\alpha=\tilde{\beta}\mu$. Take any $\nu\neq\mu$ with $|\nu|=|\mu|$ and put 
$v=S_1S_{\tilde{\beta}\mu}^*+S_2^{|\mu|}S_{\tilde{\beta}\nu}^*\;+$ (other terms). Then 
$w=vu\varphi(v^*)$ has the required form. 
\hfill$\Box$ 

\medskip
Let $u=\sum_{(\alpha,\beta)\in\J}S_\alpha S_\beta^*$ and $\ell=\ell(\J)$. Assume $\D_n^{\ell}\subseteq 
\lambda_u(\O_n)$. Then for each $(\alpha,\beta)\in\J$ and $j=1,\ldots,n$ the element 
$S_\alpha S_\beta^*S_j=P_\alpha\lambda_u(S_j)$ belongs to $\lambda_u(\O_n)$. Denote by 
${\mathcal Z}_u$ the collection of all finite products of these elements $S_\alpha S_{\tilde{\beta}}^*$ and 
their adjoints. The linear span of ${\mathcal Z}_u$ is dense in $\lambda_u(\O_n)$. Also, we denote 
by $\langle {\mathcal Z}_u \rangle$ the collection of all sums of elements from ${\mathcal Z}_u$. 

\begin{lemma}\label{finvertibility}
Let $u=\sum_{(\alpha,\beta)\in\J}S_\alpha S_\beta^*$. Denote $\ell=\ell(\J)$ and let 
$k\geq\ell$ be any integer such that there exists a $z\in{\mathcal Z}_u$, a word of length $2k-1$, with $z\in\O_n^{(1)}$. Assume that $\lambda_u(\D_n)=\D_n$ and $E(u)\neq 0$. 
Then for $\lambda_u$ to be an automorphism of $\O_n$ it suffices that $\F_n^k\subseteq\lambda_u(\O_n)$. 
If this is the case then each $S_\mu S_\nu^*$ with $\mu,\nu\in W_n^k$ belongs to 
$\langle {\mathcal Z}_u \rangle$. 
\end{lemma}
{\em Proof.}
At first we note that a word $z$, as in the statement of this lemma, exists since $E(u)\neq 0$ by assumption. 
Then observe that $\varphi^k(S_i)$ belongs to $\lambda_u(\O_n)$ for all $i=1,\ldots,n$. 
Hence $\varphi^k(\F_n)\subseteq\lambda_u(\O_n)$, 
and consequently $\lambda_u(\O_n)$ contains the entire $\F_n$. Thus $\F_n$ and $\varphi^k(S_i)$ 
are contained in $\lambda_u(\O_n)$ and we conclude that $\lambda_u(\O_n)=\O_n$. 

We have $\langle \Z_u \rangle =\lambda_u(\langle \{S_\mu S_\nu^*:\mu,\nu\in W_n\}\rangle)
\subseteq \langle \{S_\mu S_\nu^*:\mu,\nu\in W_n\}\rangle$. Now, if $\lambda_u$ is invertible then 
there exists a unitary $w\in\SS_n$ such that $\lambda_u^{-1}=\lambda_w$. Thus we have 
$$ \lambda_u^{-1}(\langle \{S_\mu S_\nu^*:\mu,\nu\in W_n\}\rangle ) 
   \subseteq \langle \{S_\mu S_\nu^*:\mu,\nu\in W_n\}\rangle $$ 
and, consequently, $\langle \Z_u \rangle = \langle \{S_\mu S_\nu^*:\mu,\nu\in W_n\}\rangle$. 
\hfill$\Box$ 

\medskip
For a while, we restrict our attention to unitaries 
$u=\sum_{(\alpha,\beta)\in\J}S_\alpha S_\beta^*$ such that $|\alpha|-|\beta|\in\{0,\pm 1\}$ 
for all $(\alpha,\beta)\in\J$. (We note that endomorphisms corresponding to such unitaries 
were studied earlier in \cite{CRS}.) 

Given a unitary $u$ as above, we may always find its presentation such that the lengths of all $\alpha$ 
coincide. Let $k$ be this common length. Then the collection of all $\alpha$ entering the presentation 
of $u$ is equal to $W_n^k$. Now, we define a {\bf new directed graph} $\Delta_u$, as follows. The set 
of vertices of $\Delta_u$ is just $W_n^k$. We put an edge from $\alpha_1$ to $\alpha_2$ 
whenever $P_{\alpha_2}
\leq P_{\tilde{\beta}_1}$. In view of our assumptions on $u$, the difference $|\alpha|-|\tilde{\beta}|$ is 
$0$, $1$ or $2$. We call this difference the degree of vertex $\alpha$ and of each edge emitted by 
$\alpha$. If $d$ is the degree of $\alpha$ then the vertex $\alpha$ emits exactly $n^d$ edges, which end 
at distinct vertices. With each edge of degree $d>0$ from $\alpha_1$ to $\alpha_2$, we associate 
a {\em label}, which is the terminal subword of length $d$ of $\alpha_2$. Edges of degree $0$ carry 
empty labels. We extend so defined labels from edges to finite directed paths on $\Delta_u$ 
by concatenation. Also, we define the {\em degree} of a path on $\Delta_u$ as the sum of the degrees 
of its edges. We denote the label of a path $x$ by $L(x)$ and its degree by $\deg(x)$. 

Now, let $\Delta_u^*$ be the set of all finite directed paths. 
In what follows, we consider pairs $(x,y)$ in $\Delta_u^*\times\Delta_u^*$ such that $x$ 
and $y$ end at the same vertex. Let $x=x'e$ and $y=y'f$, where $e$ from $\alpha_1$ to $\alpha$ 
and $f$ from $\alpha_2$ to $\alpha$ are the last edges of $x$ and $y$, respectively.  Since $e$ and $f$ end 
at the same vertex, $P_{\tilde{\beta}_1}P_{\tilde{\beta}_2}\neq0$ and thus either $\tilde{\beta}_1
\prec\tilde{\beta}_2$ or $\tilde{\beta}_2\prec\tilde{\beta}_1$. Let $\mu$ be the 
word of length $||\tilde{\beta}_1|-|\tilde{\beta}_2||$ such that  $\tilde{\beta}_1=\tilde{\beta_2}\mu$ or 
$\tilde{\beta}_2=\tilde{\beta_1}\mu$. We say that the pair $(x,y)$ is {\em balanced}  
if the following condition holds:  $L(x')=L(y')\mu$ if $\tilde{\beta}_1=\tilde{\beta_2}\mu$, and 
$L(x')\mu=L(y')$ if $\tilde{\beta}_2=\tilde{\beta_1}\mu$. Then we define the {\em total label} 
of $(x,y)$ as $L(x')$ in the former case, and $L(y')$ in the latter. Now, we define a subset $\Omega_u$ 
of the Cartesian product $\Delta_u^*\times\Delta_u^*$, as follows. A pair $(x,y)$ belongs to $\Omega_u$ 
if and only if: 
\begin{description}
\item{(i)} The paths $x$ and $y$ end at the same vertex, but they begin at distinct vertices. 
\item{(ii)} The paths $x$ and $y$ have identical degrees. 
\item{(iii)} The pair $(x,y)$ is balanced. 
\end{description}

The importance of the set $\Omega_u$ for our purposes comes from the following Proposition \ref{invcriterion}. 
Unfortunately, it is not clear to us at the present moment if its hypothesis may be algorithmically verified 
in all cases (i.e. for all applicable unitaries $u\in\SS_n$). However, in many concrete situations this can 
be done fairly easily, but preferably with the help of a computer. Thus, combined with Theorem 
\ref{diagonaltheorem}, Lemma \ref{enz} and 
Lemma \ref{finvertibility}, Proposition \ref{invcriterion} gives a criterion for deciding invertibility of 
endomorphism $\lambda_u$. 

\begin{proposition}\label{invcriterion}
Let $u=\sum_{(\alpha,\beta)\in\J}S_\alpha S_\beta^*$ be such that $|\alpha|-|\beta|\in\{0,\pm 1\}$ 
and let $\Delta_u$ be the corresponding graph. We assume that $\lambda_u(\D_n)=\D_n$. 
Then the following hold. 
\begin{enumerate}
\item Let $(x,y)\in\Omega_u$ have the total label $\gamma$ and let the paths $x,y$ begin at $\alpha$ 
and $\alpha'$, respectively. Then $\Z_u$ contains $S_{\alpha}P_\gamma S_{\alpha'}^*$. 
\item If $\alpha,\alpha'\in W_n^k$ then $S_{\alpha}S_{\alpha'}^*$ belongs to 
$\langle \Z_u \rangle$ if and only if there exists a finite collection $(x_1,y_1),\ldots,(x_m,y_m)$ in 
$\Delta_u$ with the total labels $\gamma_1,\ldots,\gamma_m$, respectively, and with all $x_j$ 
beginning at $\alpha$ and all $y_j$ beginning at $\alpha'$, such that 
$$ 1=\sum_{j=1}^m P_{\gamma_j}. $$
\end{enumerate}
\end{proposition}
{\em Proof.} 
Ad 1. Let $(\alpha,\alpha_1,\ldots,\alpha_m)$ be the consecutive vertices through which the path $x$ 
passes, and likewise let $(\alpha',\alpha_1',\ldots,\alpha_r')$ be such vertices for $y$. Then our definition 
of $\Omega_u$ ensures that 
$$ S_{\alpha}P_\gamma S_{\alpha'}^* = S_\alpha S_{\tilde{\beta}}^*S_{\alpha_1} S_{\tilde{\beta}_1}^*
\cdots S_{\alpha_m} S_{\tilde{\beta}_m}^* ( S_{\alpha'} S_{\tilde{\beta}'}^*S_{\alpha_1'} 
S_{\tilde{\beta}_1'}^*\cdots S_{\alpha_r'} S_{\tilde{\beta}_r'}^* )^*, $$
and thus $S_{\alpha}P_\gamma S_{\alpha'}^*\in\Z_u$. 

Ad 2. Suppose that $S_{\alpha}S_{\alpha'}^*\in\langle \Z_u \rangle$, and let $S_{\alpha}S_{\alpha'}^* = 
\sum_{j=1}^m S_{\mu_j}S_{\nu_j}^*$, with each $S_{\mu_j}S_{\nu_j}^*$ in $\Z_u$. Since there are 
no cancellations among words, each $S_{\mu_j}S_{\nu_j}^*$ must be of the form $S_\alpha P_{\gamma_j} 
S_{\alpha'}^*$ for some $\gamma_j\in W_n$. Now, it is not difficult to verify that an element of $\Z_u$ 
has this form if and only if there exists a pair $(x_j,y_j)$ in $\Omega_u$ with the total label $\gamma_j$ and 
such that $x_j$ and $y_j$ begin at $\alpha$ and $\alpha'$, respectively.  

The reverse implication is an immediate consequence of part 1 of this proposition. 
\hfill$\Box$

\medskip
We end this section with some examples of invertible endomorphisms $\lambda_u$, $u\in\SS_n\setminus\P_n$.

\begin{example}\label{realexam}
\rm Let $\mu,\nu$ be two words such that $\tilde{\nu}=j_1\cdots j_r\tilde{\mu}$ with $j_k\in W_n^1$ 
and $j_k\not\in\{\mu_1,\nu_1\}$ for all $k=1,\ldots,r$. Let 
$$ u = S_\nu S_\mu^* + S_\mu S_\nu^* + 1-P_\nu-P_\mu. $$
Suppose that $\lambda_u(\D_n)=\D_n$. We claim that then $\lambda_u$ is automatically invertible. 
Indeed, it suffices to check that $S_\nu S_\mu^*\in\lambda_u(\O_n)$. But we have 
$S_\nu S_\mu^* = S_\nu S_{\tilde{\mu}}^* S_{\tilde{\mu}} S_{\tilde{\nu}}^* S_{\tilde{\nu}} S_\mu^*$. 
Now, $S_\nu S_{\tilde{\mu}}^* = P_\nu\lambda_u(S_{\mu_1})$ and $S_{\tilde{\nu}} S_\mu^*= 
\lambda_u(S_{\nu_1}^*)P_\mu$ are both in $\lambda_u(\O_n)$. Also, $S_{\tilde{\nu}}
S_{\tilde{\mu}}^* = S_{\tilde{\nu}}S_{\tilde{\nu}}^*S_{j_1}\cdots S_{j_r} = P_{\tilde{\nu}} 
\lambda_u(S_{j_1})\cdots\lambda_u(S_{j_r})$ and hence $S_{\tilde{\mu}} S_{\tilde{\nu}}^*\in
\lambda_u(\O_n)$. Consequently, $\lambda_u$ is invertible, as claimed. 
\end{example}

\begin{example}\label{realexam2}
\rm Let $\alpha_1,\alpha_2,\alpha_3$ be such that $\{P_{\alpha_j}\}$ are mutually orthogonal and each $\alpha_j$ 
begins with the same letter $i$. Furthermore, suppose that $\tilde{\alpha}_j=\gamma_j\mu$ for some 
$\gamma_j$ which do not contain the letter $i$. Let 
$$ u=S_{\alpha_1}S_{\alpha_2}^* + S_{\alpha_2}S_{\alpha_3}^* + S_{\alpha_3}S_{\alpha_1}^* + 
1- P_{\alpha_1}-P_{\alpha_2}-P_{\alpha_3}. $$
If $\lambda_u(\D_n)=\D_n$ then automatically $\lambda_u
\in\aut(\O_n)$. Indeed, we have $S_k=\lambda_u(S_k)$ for all $k\neq i$ and thus 
$S_{\alpha_1}S_{\alpha_2}^*=P_{\alpha_1}\lambda(S_i)S_{\gamma_2}S_{\gamma_3}^*\lambda(S_i^*)
P_{\alpha_2}$ belongs to $\lambda_u(\O_n)$. Similarly, $S_{\alpha_2}S_{\alpha_3}^*$ and $S_{\alpha_3}S_{\alpha_1}^*$ are in $\lambda_u(\O_n)$ as well. Thus $u\in\lambda_u(\O_n)$ and 
$\lambda_u$ is invertible. A concrete example in $\O_2$ is obtained by putting
$$ w=S_{11}S_{121}^* + S_{121}S_{1221}^* + S_{1221}S_{11}^* + P_{1222} + P_2, $$ 
and then indeed $\lambda_w$ is an automorphism of $\O_2$. 
\end{example}


\section{The outer Weyl group}

In this section, we consider the question whether an endomorphism corresponding to a unitary in $\SS_n$ 
may or may not be equivalent (via an inner automorphism) to one corresponding to a unitary in the core 
UHF-subalgebra $\F_n$. 

\begin{proposition}\label{nonfnconjugate}
There exist unitaries $u\in\SS_n$ such that $\lambda_u\not\in\aut(\O_n)\lambda(\F_n)$. 
\end{proposition}
{\em Proof.}
At first we observe that if $w\in\U(\F_n)$ and $Q\neq 0$ is a projection in $\O_n$ then the space $\lambda_w(\D_n)Q$ is infinite dimensional. Indeed, since $E(Q)$ is a non-zero, positive element of $\F_n$, 
there is a non-zero projection $q\in\F_n$ and a scalar $t>0$ such that $tq\leq E(Q)$. There exists 
a sequence of indices $j_k\in W_n^1$ such that if $\alpha_k\in W_n^k$ are defined recursively as 
$\alpha_1=j_1$, $\alpha_{k+1}=(\alpha_k,j_{k+1})$ then $\lambda_w(P_{\alpha_k})q\neq 0$ for 
all $k$. The sequence $\{q\lambda_w(P_{\alpha_k})q\}$ never stabilizes. Indeed, if  $q\lambda_w(P_{\alpha_{k+m}})q=q\lambda_w(P_{\alpha_k})q$ for all $m$ then 
$$ 0 \neq \tau(q\lambda_w(P_{\alpha_k})q)=\tau(q\lambda_w(P_{\alpha_{k+m}})q) = 
   \tau(\lambda_w(P_{\alpha_{k+m}})q\lambda_w(P_{\alpha_{k+m}})) \leq \tau(P_{\alpha_{k+m}}) 
   \underset{m\to\infty}{\longrightarrow} 0, $$ 
a contradiction. The inequality above holds since $w$ being in $\U(\F_n)$ the corresponding 
endomorphism $\lambda_w$ is $\tau$-preserving. Thus, there is a strictly decreasing, infinite sequence of 
projections $f_1>f_2>\ldots$ in $\lambda_w(\D_n)$ such that $qf_kq>qf_{k+1}q$ for all $k$. 
Thus $(f_k-f_{k+1})q\neq0$ for all $k$, and hence 
$$ 0 \neq (f_k-f_{k+1})tq(f_k-f_{k+1}) \leq (f_k-f_{k+1})E(Q)(f_k-f_{k+1}). $$ 
Thus  $(f_k-f_{k+1})Q\neq0$ and, consequently, $\{(f_k-f_{k+1})Q\}$ is an infinite sequence 
of linearly independent elements of $\lambda_w(\D_n)Q$, since these are non-zero operators 
with mutually orthogonal ranges. 

Now, the same conclusion as above holds if $\lambda_w$ is replaced by $\psi\lambda_w$ for 
some automorphism $\psi\in\aut(\O_n)$, since the dimension of $(\psi\lambda_w)(\D_n)Q$ is 
the same as that of $\lambda_w(\D_n)\psi^{-1}(Q)$. Thus, the conclusion of the proposition follows 
from Example \ref{examples123} (a), where a unitary $u\in\SS_2$ is exhibited such that $\lambda_u(\D_2)
P_2$ is one-dimensional.  
\hfill$\Box$

\medskip
Of course, the method of Proposition \ref{nonfnconjugate} cannot give any information about 
automorphisms. We treat the automorphism case in Theorem \ref{outermain}, below. To the best of our 
knowledge, the automorphism entering its proof is the first known example of an automorphism 
of $\O_n$ not inner related to an automorphism induced by a unitary from 
the core UHF subalgebra $\F_n$. 

\begin{theorem}\label{outermain}
There exist automorphisms $\lambda_u$, $u\in\SS_n$, of $\O_n$ such that for all $w\in\U(\O_n)$ and  
$v\in\U(\F_n)$ we have $\lambda_u \neq \Ad(w)\lambda_v$. 
\end{theorem}
{\em Proof.}
At first we consider the following self-adjoint element of $\SS_2$ (c.f. Example \ref{realexam}):  
\begin{equation}\label{superu}
u=S_{11}S_{121}^* + S_{121}S_{11}^* + P_{122} + P_2. 
\end{equation}
One easily checks that 
$$ \lambda_u(S_1) = S_1(S_1S_{21}^* + S_{21}S_1^* +P_{22}) \;\;\; \text{and} \;\;\;
\lambda_u(S_2) = S_2. $$
This yields $\lambda_u^2=\id$. Suppose by contradiction that there exist $w\in\U(\O_2)$ 
and $v\in\U(\F_2)$ such that $\lambda_u=\Ad(w)\lambda_v$. Then we have $\D_2=\lambda_u(\D_2)=
\Ad(w)\lambda_v(\D_2)$ and $\lambda_v(\D_2)\subseteq\F_2$, thus $w^*\D_2 w\subseteq\F_2$. 
Hence for all $d\in\D_2$ and $z\in\operatorname{U(1)}$ we have $\gamma_z(w^*dw)=w^*dw$. 
Therefore $\gamma_z(w)w^*\in\D_2'\cap\O_2=\D_2$. Thus for each $z\in\operatorname{U(1)}$ there  
exists a unitary $d_z\in\U(\D_2)$ such that $\gamma_z(w)=d_z w$. Now we calculate
$$ 
\gamma_z(u) = \gamma_z(wv\varphi(w^*))  = d_z wv\varphi(w^*)\varphi(d_z^*) =d_z u\varphi(d_z^*). 
$$
This means that if $u=\sum S_\alpha S_\beta^*$ then  
$$ z^{|\alpha|-|\beta|}S_\alpha S_{\tilde{\beta}}^*d_z=d_z S_\alpha S_{\tilde{\beta}}^* $$ 
for all $(\alpha,\beta)$ and $z\in\operatorname{U(1)}$. Taking $\alpha=\beta=2$ this yields 
$S_2 d_z =  d_z S_2$. By \cite{MT}, $d_z$ is a scalar, and this in turn implies that 
$\gamma_z(u)=u$ for all $z$, a contradiction.

Now, if $n\geq2$ is arbitrary, then we consider $\tilde{u}=S_{11}S_{121}^* + S_{121}S_{11}^* + 
P_{122} + 1-P_1$, and the same argument as above applies. 
\hfill$\Box$

\medskip
As immediate consequences of Theorem \ref{outermain}, we obtain the following two corollaries. 

\begin{corollary}\label{fullrestricted}
The restricted outer Weyl group of $\O_n$ is a proper subgroup of the outer Weyl group of $\O_n$. 
\end{corollary}

As shown in \cite{CHS1}, the restricted outer Weyl group of $\O_n$ is residually finite and nonamenable. 
Thus the outer Weyl group is nonamenable as well, but we do not know if it is residually finite. 

\begin{corollary}\label{innerouterpairconjugate}
There exist unital subalgebras $\A$ of $\O_n$ isomorphic to the UHF algebra of type $\{n^\infty\}$  
such that $\F_n$ and $\A$ are conjugate inside $\O_n$ (by an automorphism of $\O_n$) 
but not inner conjugate. 
\end{corollary}


\section{The action on the shift space}

Equality (\ref{adaction}) easily implies that for all $d\in\D_n$ and all $k>\ell'(\J)$ we have 
\begin{equation}\label{adshift}
\Ad(u)(\varphi^k(d)) = \sum_{(\alpha,\beta)\in\J}\varphi^{k+|\alpha|-|\beta|}(d)P_\alpha. 
\end{equation}

Consider a map $f:\D_n\to\D_n$. We say that $f$ eventually preseves standard projections if there 
exists an integer $m\in\Nb$ such that for each $\alpha\in W_n$, $|\alpha|\geq m$, the image $f(P_\alpha)$ 
is a standard projection. If $u\in\SS_n$ then $\Ad(u)$ eventually preserves standard projections. 

\begin{proposition}\label{adactionchar}
If $f\in\aut(\D_n)$ then there exists a unitary $u\in\SS_n$ such that $f=\Ad(u)|_{\D_n}$ if and only if; 
\begin{description}
\item{(i)} $f$ eventually preserves standard projections, and 
\item{(ii)} there exist projections $P_i,Q_i$, $i=1,\ldots,r$, in $\D_n$ and non-negative integers 
$k_i,m_i$, $i=1,\ldots,r$, such that $\sum_{i=1}^r P_i=1=\sum_{i=1}^r Q_i$ and  
$$ f\circ M_{P_i}\circ\varphi^{k_i} = M_{Q_i}\circ\varphi^{m_i}, \;\;\;\;\; i=1,\ldots,r. $$
\end{description}
\end{proposition}
{\em Proof.} 
Let $f\in\aut(\D_n)$ satisfy conditions (i) and (ii) of the proposition. For a given $i\in\{1,\ldots,r\}$, we  
note that for any subprojection $p$ of $P_i$ we have $ f\circ M_p\circ\varphi^{k_i} = M_{f(p)}
\circ\varphi^{m_i}$. Subdividing $P_i$ into a sum of standard projections and using condition (i), 
we can assume in condition (ii) that all projections $P_i,Q_i$ are standard, say $P_i=P_{\beta_i}$ 
and $Q_i=P_{\alpha_i}$. Define $u=\sum_{i=1}^r S_{\alpha_i} S^*_{\beta_i}$, a unitary 
element of $\SS_n$. Then we have 
$$ (\Ad(u^*)\circ f)\circ M_{P_{\beta_i}}\circ\varphi^{k_i+|\alpha_i|+h} 
= M_{P_{\beta_i}}\circ\varphi^{m_i+|\beta|+h}, \;\;\;\;\; i=1,\ldots,r, $$
for all sufficiently large $h\in\Nb$. We claim that $k_i+|\alpha_i| = m_i+|\beta_i|$ 
for each $i$. Indeed, fix an $i$ and suppose that $k_i+|\alpha_i| \geq m_i+|\beta_i|$ 
(otherwise consider $ (\Ad(u^*)\circ f)^{-1}$ instead). Then we have 
$$ (\Ad(u^*)\circ f)(\varphi^{k_i+|\alpha_i|-m_i-|\beta_i|}(y)P_{\beta_i})=yP_{\beta_i}, \;\;\; 
\forall y\in\varphi^{m_i+|\beta_i|+h}(\D_n), $$
for all sufficiently large $h\in\Nb$. Fix such an $h$ and let $r\geq m_i+|\beta_i|+h$ be such that 
$$  (\Ad(u^*)\circ f)(\varphi^{k_i+|\alpha_i|-m_i-|\beta_i|}(\D_n^{m_i+|\beta_i|+h})P_{\beta_i}) 
\subseteq \D_n^r P_{\beta_i} $$
and $P_{\beta_i}\in\D_n^r$. Then we have 
$$ \begin{aligned}
 & (\Ad(u^*)\circ f)(\varphi^{k_i+|\alpha_i|-m_i-|\beta_i|}(\D_n^r)P_{\beta_i}) \\ 
 & \phantom{xxxxx} = (\Ad(u^*)\circ f)(\varphi^{k_i+|\alpha_i|-m_i-|\beta_i|}(\D_n^{m_i+|\beta_i|+h})
\varphi^{m_i+|\beta_i|+h}(\D_n^{r-m_i-|\beta_i|-h})P_{\beta_i}) \\
 & \phantom{xxxxx}\subseteq \D_n^r P_{\beta_i} \varphi^{m_i+|\beta_i|+
h}(\D_n^{r-m_i-|\beta_i|-h})P_{\beta_i} \\
 & \phantom{xxxxx} \subseteq \D_n^r P_{\beta_i}.
\end{aligned} $$
Since $\Ad(u^*)\circ f$ is injective and the dimension of $\varphi^{k_i+|\alpha_i|-m_i-|\beta_i|}(\D_n^r)
P_{\beta_i}$ is not smaller than the dimension of $\D_n^r P_{\beta_i}$, it follows that these two 
dimensions are identical, and this can only happen when $k_i+|\alpha_i|-m_i-|\beta_i|=0$. 
Consequently, 
$$ (\Ad(u^*)\circ f)\circ M_{P_{\beta_i}}\circ\varphi^{h} 
= M_{P_{\beta_i}}\circ\varphi^{h}, \;\;\;\;\; i=1,\ldots,r, $$
for all sufficiently large $h\in\Nb$. Summing over $i$ we get 
$$ (\Ad(u^*)\circ f)\circ\varphi^{h} = \varphi^{h} $$
for all sufficiently large $h\in\Nb$. Therefore $\Ad(u^*)\circ f = \Ad(w)|_{\D_n}$ for some $w\in\P_n$, 
by \cite[Lemma 3.2]{CHS1}. Hence $f=\Ad(uw)|_{\D_n}$ and $uw\in\SS_n$. 
This proves one direction. The opposite direction is clear. Indeed, let $u\in\SS_n$ be such that 
$u=\sum_{(\alpha,\beta)\in\J}S_\alpha S_\beta^*$. Then condition (i) holds, as noted just above 
this proposition. One easily checks that condition (ii) holds with projections $P_\beta$ and $P_\alpha$ 
instead of $P_i$ and $Q_i$, respectively, and with $|\beta|$ and $|\alpha|$ instead of $k_i$ 
and $m_i$, respectively. 
\hfill$\Box$

\medskip
Given $u\in\SS_n$ and considering the homeomorphism $\Ad(u)_*$ of the spectrum $X_n$ of $\D_n$, 
we see that the set of fixed points has a very simple structure, as the following Proposition 
\ref{innerfixedpoints} shows. 

\begin{proposition}\label{innerfixedpoints}
For $u\in\SS_n$, the set of fixed points in $X_n$ for the homeomorphism $\Ad(u)_*$ consists of the 
union of a clopen set and a finite set. Furthermore, each of the isolated fixed points is either a local attractor 
or a local repeller. 
\end{proposition}
{\em Proof.}
Let $u=\sum_{(\alpha,\beta)\in\J}S_\alpha S_\beta^*$. It is clear that $\Ad(u)_*$ admits fixed points 
in $X_n$ if and only if there exists $(\alpha,\beta)\in\J$ such that either $\alpha\prec\beta$ or 
$\beta\prec\alpha$. Thus we arrive at one of the following three cases. (1) If $\alpha=\beta$ then 
the clopen set $\{x\in X_n : \beta\prec x\}$ is fixed by  $\Ad(u)_*$. (2) If $\alpha=\beta\mu$, 
$\mu\neq\emptyset$, then $x=\beta\mu\mu\ldots$ is a fixed point and a local attractor. (3) If 
$\beta=\alpha\mu$, $\mu\neq\emptyset$, then $x=\alpha\mu\mu\ldots$ is a fixed point and a local 
repeller. 
\hfill$\Box$

\medskip
In contrast to Proposition \ref{innerfixedpoints} above, the set of fixed points in $X_n$ corresponding to 
an outer automorphism $\lambda_u$, $u\in\SS_n$, may have a much more complicated structure, as the 
following example demonstrates. 

\begin{example}\label{superufixedpoints}
\rm Let $u$ be the unitary in $\SS_2$ defined by formula (\ref{superu}). It is not difficult to verify 
that the corresponding homeomorphism $(\lambda_u)_*$ of $X_2$ fixes an $x\in X_2$ if and 
only if $x$ does not contain substrings $(11)$ and $(121)$. These fixed points form a compact, 
nowhere dense subset $K$ of $X_2$, in which there are no isolated points. Thus $K$ itself is 
homeomorphic to the Cantor set and closed under the action of the one-sided shift $\varphi_*$. 
\end{example}


\medskip\noindent
Roberto Conti \\
Dipartimenti di Scienze \\
Universit{\`a} di Chieti-Pescara `G. D'Annunzio' \\
Viale Pindaro 42, I--65127 Pescara, Italy 

\vspace{2mm}\noindent
present address: \\
Dipartimento di Scienze di Base e Applicate per l'Ingegneria \\
Sezione di Matematica \\
Sapienza Universit\`a di Roma \\
Via A. Scarpa 16 \\
00161 Roma, Italy \\
E-mail: roberto.conti@sbai.uniroma1.it \\

\smallskip\noindent
Jeong Hee Hong \\
Department of Data Information \\
Korea Maritime University \\
Busan 606--791, South Korea \\
E-mail: hongjh@hhu.ac.kr \\

\smallskip \noindent
Wojciech Szyma{\'n}ski\\
Department of Mathematics and Computer Science \\
The University of Southern Denmark \\
Campusvej 55, DK-5230 Odense M, Denmark \\
E-mail: szymanski@imada.sdu.dk

\end{document}